\documentclass[leqno,11pt]{article} 
 
\usepackage{amsmath,amssymb,amsbsy,amsfonts,amsthm,latexsym,
            amsopn,amstext,amsxtra,euscript,amscd,amsthm}
\usepackage{amsmath} 
\usepackage{amsthm} 
\usepackage{latexsym}
\usepackage{amssymb} 
\usepackage{xcolor}
 
\usepackage[applemac]{inputenc}
\usepackage[pdftex]{graphicx}
\usepackage
{hyperref} 

\usepackage{amsbsy}

\usepackage{fancyhdr}
 
\def\boxit#1#2{\setbox1=\hbox{\kern#1{#2}\kern#1}%
\dimen1=\ht1 \advance\dimen1 by #1 \dimen2=\dp1 \advance\dimen2 by #1
\setbox1=\hbox{\vrule height\dimen1 depth\dimen2\box1\vrule}%
\setbox1=\vbox{\hrule\box1\hrule}%
\advance\dimen1 by .4pt \ht1=\dimen1
\advance\dimen2 by .4pt \dp1=\dimen2 \box1\relax}

\def\hfl#1#2{\smash{\mathop{\hbox to 12 mm{\rightarrowfill}}
\limits^{\scriptstyle#1}_{\scriptstyle#2}}}

 \def\hlfl#1#2{\smash{\mathop{\hbox to 12 mm{\leftarrowfill}}
\limits^{\scriptstyle#1}_{\scriptstyle#2}}}

\def\phfl#1#2{\smash{\mathop{\hbox to 8 mm{\rightarrowfill}}
\limits^{\scriptstyle#1}_{\scriptstyle#2}}}

 \def\phlfl#1#2{\smash{\mathop{\hbox to 8 mm{\leftarrowfill}}
\limits^{\scriptstyle#1}_{\scriptstyle#2}}}

\def\cqfd{\unskip\kern 6pt\penalty 500
\raise -2pt\hbox{\vrule\vbox to 10pt{\hrule width 4pt
\vfill\hrule}\vrule}\par}

\def\house#1{\setbox1=\hbox{$\,#1\,$}%
\dimen1=\ht1 \advance\dimen1 by 2pt \dimen2=\dp1 \advance\dimen2 by 2pt
\setbox1=\hbox{\vrule height\dimen1 depth\dimen2\box1\vrule}%
\setbox1=\vbox{\hrule\box1}%
\advance\dimen1 by .4pt \ht1=\dimen1
\advance\dimen2 by .4pt \dp1=\dimen2 \box1\relax}

\def\virgule{\raise 2pt \hbox{$,$}}
 
\def\C{\mathbb{C}}

\def\Q{\mathbb{Q}}
\def\R{\mathbb{R}}
\def\Z{\mathbb{Z}}

\def\rmh{{\mathrm{h}}}
\def\rmH{{\mathrm{H}}}
\def\rmM{{\mathrm{M}}}
\def\rmN{{\mathrm{N}}}

\def\rmM{{\mathrm{M}}}

\def\og{\leavevmode\raise.3ex\hbox{$\scriptscriptstyle
\langle\!\langle$}}

\def\fg{\leavevmode\raise.3ex\hbox{$\scriptscriptstyle
\,\rangle\!\rangle$}}

\def\Qbar{\overline{\Q}}
\def\betatilde{\tilde{\beta}}

\newtheorem{theorem}{\indent  Theorem }
\newtheorem{proposition}[theorem]{\indent  Proposition}
\newtheorem{lemma}[theorem]{\indent  Lemma}
\newtheorem{corollary}[theorem]{\indent  Corollary}

\newtheorem*{namedtheorem}{\\indent  theoremname}

\newcounter{compteurkappa} 

\def\Newcst#1{
\refstepcounter{compteurkappa}
\kappa_{ 
\arabic{compteurkappa}}
\label{#1}
}

\def\cst#1{\kappa_{\ref{#1}}}

\begin{document}
 

 \begin{center}
 {\bf \Large
Families of Thue equations associated 
\\
with a 
rank one subgroup of the 
\\
unit group of a number field
\\

}

 \bigskip
{\large by}

 \bigskip
{ \it \Large
Claude Levesque } {\large and} {\it \Large
 Michel Waldschmidt
 }
 \end{center}\bigskip\medskip

{\sc Abstract.}
Twisting a binary form $F_0(X,Y)\in\Z[X,Y]$ of degree $d\ge 3$ by powers $\upsilon^a$ ($a\in\Z$) of an algebraic unit $\upsilon$ gives rise to a binary form $F_a(X,Y)\in\Z[X,Y]$. More precisely, when $K$ is a number field of degree $d$, $\sigma_1,\sigma_2,\dots,\sigma_d$ the embeddings of $K$ into $\C$, $\alpha$ a nonzero element in $K$, $a_0\in\Z$, $a_0>0$ and \\
\centerline{$
 F_0(X,Y)=a_0\displaystyle\prod_{i=1}^d (X-\sigma_i(\alpha) Y),
$}
then for $a\in\Z$ we set\\
\centerline{$
 F_a(X,Y)=\displaystyle a_0\prod_{i=1}^d (X-\sigma_i(\alpha\upsilon^a) Y).
 $}
 Given $m\ge 0$, our main result is an effective upper bound for the solutions $(x,y,a)\in\Z^3$ of the Diophantine inequalities\\
\vskip -3mm
\centerline{$
0<|F_a(x,y)|\le m
$}
\noindent for which $xy\not=0$ and $\Q(\alpha \upsilon^a)=K$.
Our estimate involves an effectively computable constant depending only on $d$; it is explicit in terms of $m$, in terms of the heights of $F_0$ and of $\upsilon$, and in terms of the regulator of the number field $K$.
\\


\section{Introduction and the main results}\label{S:introduction}
 
Let $d\ge 3$ be a given integer. We denote by $\kappa_1,\kappa_2,\dots$ positive effectively computable constants which depend only on $d$.

Let $K$ be a number field of degree $d$. Denote by $\sigma_1,\sigma_2,\dots,\sigma_d$ the embeddings of $K$ into $\C$ and by $R$ the regulator of $K$. Let $\alpha\in K$, $\alpha\not=0$, and let $a_0\in\Z$, $a_0>0$, be such that the coefficients of the polynomial
$$
f_0(X)=a_0\prod_{i=1}^d \bigl(X-\sigma_i(\alpha)\bigr)
$$
are in $\Z$. Let $\upsilon$ be a unit in $K$, not a root of unity. For $a\in\Z$, define the polynomial $f_a(X)$ in $\Z[X]$ and the binary form $F_a(X,Y)$ in $\Z[X,Y]$ by
$$
f_a(X)=a_0\prod_{i=1}^d \bigl(X-\sigma_i(\alpha \upsilon^a)\bigr)
$$
and
$$
F_a(X,Y)=Y^df_a(X/Y)=a_0\prod_{i=1}^d \bigl(X-\sigma_i(\alpha \upsilon^a)Y\bigr).
$$
Define
$$
\lambda_0=a_0\prod_{i=1}^d\max\{1,|\sigma_i(\alpha)|\}
\quad\hbox{and}\quad
\lambda=\prod_{i=1}^d\max\{1,|\sigma_i(\upsilon)|\}.
$$

Let $m\in\Z$, $m>0$. We consider the family of Diophantine inequalities
\begin{equation}\label{Equation:ThueInequality}
0<|F_a(x,y)|\le m,
\end{equation}
where the unknowns $(x,y,a)$ take their values in the set of elements in $\Z^3$ such that $xy\not=0$ and $\Q(\alpha \upsilon^a)=K$. It follows from the results in \cite{LW1} that the set of solutions is finite. However, the proof in \cite{LW1} relies on Schmidt's subspace theorem, which is not effective. Here we give an effective upper bound for $\max\{|x|,|y|,|a|\}$ in terms of $m$, $R$, $\lambda_0$ and $\lambda$, by using lower bounds for linear forms in logarithms. 

For $x\in\R$, $x>0$, we stand to the notation $\log^\star x$ for $\max\{1,\log x\}$. 

Here is our main result. 

\begin{theorem}\label{Theorem:mainGeneralWeak} 
There exists an effectively computable constant $\Newcst{kappa:mainTh}>0$, depending only on $d$, 
such that any solution $(x,y,a)\in\Z^3$ of $(\ref{Equation:ThueInequality})$,
which verifies $xy \neq 0$ and $\Q(\alpha \upsilon^a)=K$, satisfies
$$
 |a| \le \cst{kappa:mainTh} \lambda^{d^2(d+2)/2} (R+\log m+\log\lambda_0)R\log^\star R .
$$ 
\end{theorem}

Under the assumptions of Theorem $\ref{Theorem:mainGeneralWeak}$, with the help of the upper bound
$$
 \rmH(F_a)\le 2^d\lambda_0\, \lambda^{|a|}
 $$
 for the height of the form $F_a$, it follows from the bound (3.2) in \cite[Theorem 3]{BG2debut} (see also \cite[Th.~9.6.2]{EG}) that
 $$
 \log \max \{ |x|, |y| \}
\le \kappa
\bigl(R+ \log^\star m+
|a|\log \lambda + \log\lambda_0
 \bigr)R (\log^\star R) 
$$
with 
$$
\kappa=3^{r+27}(r+1)^{7r+19}d^{2d+6r+15}.
$$
Combining this upper bound with our Theorem $\ref{Theorem:mainGeneralWeak}$ provides an effective upper bound 
for $\max \{ |x|, |y|,|a| \}$.

For $i=1,\dots,d$, set $\upsilon_i=\sigma_i(\upsilon)$ and assume
$$
|\upsilon_1|\le |\upsilon_2|\le \cdots \le |\upsilon_d|.
$$ 
Our proof actually gives a much stronger estimate for $|a|$, see Theorem $\ref{Theorem:majorationdea}$, which involves some extra parameter $\mu>1$ defined by 
$$
\mu=\begin{cases} 
 \lambda
& \hbox{if 
$ |\upsilon_1|= |\upsilon_{d-1}|$ or $ |\upsilon_2|= |\upsilon_d| $,}
\\[1mm]
\displaystyle \min \left\{ \frac{ |\upsilon_{d-1}|}{ |\upsilon_1|}, \frac{ |\upsilon_d|}{ |\upsilon_2|}\right\}
& \hbox{if 
$|\upsilon_1|<|\upsilon_2|= |\upsilon_{d-1}| <|\upsilon_d|$,}
\\[3mm]
\displaystyle 
\frac{ |\upsilon_{d-1}|}{ |\upsilon_2|} & \hbox{if 
$ |\upsilon_2|< |\upsilon_{d-1}| $.}\\
\end{cases}
$$
Notice that the condition 
$ |\upsilon_1|= |\upsilon_{d-1}|$ means $ |\upsilon_1|=|\upsilon_2|=\cdots= |\upsilon_{d-1}|$ and that the condition $ |\upsilon_2|= |\upsilon_d|$ means $ |\upsilon_2|=|\upsilon_3|=\cdots= |\upsilon_d|$; using Lemma $\ref{Lemma:deuxconjuguesreels}$, we deduce that each of these two conditions implies that $d$ is odd, hence that the field $K$ is almost totally imaginary (namely, with a single real embedding) -- compare with \cite{LW-Balu}. 

\begin{theorem}\label{Theorem:majorationdea}
There exists a positive effectively computable constant $\Newcst{kappa:maja1}$, depending only on $d$, with the following property. 
Let $(x,y,a)\in \Z^3$ satisfy 
$$
xy\not=0, \quad[\Q(\alpha \upsilon^a):\Q]=d
 \quad\hbox{and}\quad 
0<|F_{a}(x,y)|\le m.
$$
Then
\begin{equation}\label{Equation:bornepouraavecmu}
 |a| 
\le \cst{kappa:maja1} \frac{\log \lambda}{\log\mu}(R+\log m+\log \lambda_0+\log \lambda) 
R \log \left( R \frac{(\log\lambda)^2}{\log \mu}\right).
\end{equation}
\end{theorem}

On the one hand, using Lemma $\ref{Lemma:minorationalpha2suralphad}$ (\S$\ref{SS:Conjugates}$), we will prove in \S$\ref{S:ProofMainTheorem}$ that 
$$
\log\mu\ge \Newcst{kappa:maja2} \lambda^{-d^2(d+2)/2}(\log\lambda)^2,
$$
which will enable us to deduce Theorem $\ref{Theorem:mainGeneralWeak}$ from Theorem $\ref{Theorem:majorationdea}$. On the other hand, thanks to $(\ref{Equation:hauteurunite}$), we have $\mu\le \lambda^2$. In general, we expect $\mu$ to be as large as $\lambda^{\Newcst{kappa:maja3} }$ (which is therefore the maximum possible), in which case the conclusion of Theorem $\ref{Theorem:majorationdea}$ becomes 
\begin{equation}\label{Equation:strongEstimate}
 |a| 
\le \cst{kappa:maja4} (R+\log m+\log\lambda_0+\log\lambda)R(\log R+\log^\star\log^\star\lambda)
\end{equation}
with a positive effective constant $ \Newcst{kappa:maja4}$ depending only on $d$.
In \S$\ref{S:example}$, we give a few examples where this last bound is valid. 

In Theorem $\ref{Theorem:mainGeneralWeak}$, the hypothesis that $\upsilon$ is not a root of unity cannot be omitted. Here is an example with $\alpha=a_0=m=1$. Let $\Phi_n(X)$ be the cyclotomic polynomial of index $n$ and degree $\varphi(n)$ (Euler totient function). Let $\zeta_n$ be a primitive $n$--th root of unity. Set $f_0=\Phi_n$ and $u=\zeta_n$. 
For $a\in\Z$ with $\gcd(a,n)=1$, the irreducible polynomial $f_a$ of $\zeta_n^a$ is nothing else than $f_0$. Hence, if the equation 
$$
F_0(x,y)=\pm 1
$$
has a solution $(x,y)\in\Z^2$ with $xy\not=0$, then for infinitely many $a\in\Z$ the twisted Thue equation $F_a(x,y)=\pm 1$ has also the solution $(x,y)$, since $F_a=F_0$. For instance, when $n=12$, we have $\Phi_{12}(X)=X^4-X^2+1$ and
the equation 
$$ 
x^4-x^2y^2+y^4=1
$$
has the solutions $(1,1)$, $(-1,1)$, $(1,-1)$, $(-1,-1)$. 

The main result of \cite{LWvdP}, which deals only with non totally real cubic equations, is a special case of Theorem $\ref{Theorem:majorationdea}$; the ``constants'' in \cite{LWvdP} depend on $\alpha$ and $\upsilon$, while here they depend only on $d$. The main result of \cite{LW2} deals with Thue equations twisted by a set of units which is not supposed to be a group of rank $1$, but it involves an assumption (namely that at least two of the conjugates of $\upsilon$ have a modulus as large as a positive power of $\house{\upsilon}$) which we do not need here. Our Theorem $\ref{Theorem:majorationdea}$ also improves the main result of \cite{LW-Ram}: we remove the assumption that the unit is totally real (besides, the result of \cite{LW-Ram} is not explicit in terms of the heights and regulator). We also notice that the part ${\mathrm{(iii)}}$ of Theorem 1.1 of \cite{LW-Thue} follows from our Theorem $\ref{Theorem:majorationdea}$. The main result of \cite{LW-Balu} does not assume that the twists are done by a group of units of rank $1$, but it needs a strong assumption which does not occur here, namely that the field $K$ has at most one real embedding. 

We conclude this \S \ref{S:introduction} with some more definitions and properties.

When $f$ is a polynomial in one variable of degree $d$ with coefficients in $\Z$ and leading coefficient $c_0>0$, the (usual) height $\rmH(f)$ of $f$ is the maximum of the absolute values of the coefficients of $f$, while the Mahler measure of $f$ is 
$$
\rmM(f)=c_0\prod_{i=1}^d \max\{1,|\gamma_i|\},
$$
where $\gamma_1,\gamma_2,\dots,\gamma_d$ are the roots of $f$ in $\C$.

Let us recall\footnote{Our $\rmh$ is the same as in \cite{EG}, it corresponds to the logarithm of the $\rmh$ in \cite{BG2debut}.}
 that the logarithmic height $\rmh(\gamma)$ of an algebraic number $\gamma$ of degree $d$ is $\frac{1}{d}\log\rmM(\gamma)$ where $\rmM(\gamma)$ is the Mahler measure of the irreducible polynomial of $\gamma$. 
We have 
\begin{equation}\label{equation:heights}
\rm M(f)\le \sqrt{d+1} \,\rm H(f)\;
\; \hbox{ and }\;\;
\rmH(f)\le 2^d 
\rmM(f) 
\end{equation} 
(see \cite{GL326}, Annex to Chapter 3, {\it Inequalities Between Different Heights of a Polynomial,} pp. 113--114; see also \cite[\S 1.9]{EG}). The second upper bound in $(\ref{equation:heights})$ could be replaced by the sharper one
$$
\rmH(f)\le 
\binom{d }{\lfloor d/2\rfloor}
\rmM(f),
$$
but we will not need it. 

Let $\upsilon$ be a unit of degree $d$ and conjugates $\upsilon_1,\dots,\upsilon_d$ with 
$$
 |\upsilon_1|\le |\upsilon_2|\le \cdots \le |\upsilon_d|, 
 $$
 so that $\house{\upsilon}=|\upsilon_d|$. 
 Let $\lambda=\rmM(\upsilon)$ and let $s$ be an index in $\{1,\dots,d-1\}$ such that
 $$
 |\upsilon_1|\le |\upsilon_2|\le \cdots \le |\upsilon_s|\le 1\le |\upsilon_{s+1}|\le \cdots \le |\upsilon_d|.
 $$
 We have 
 $$
\lambda= \rmM(\upsilon)=|\upsilon_{s+1} \cdots \upsilon_d|\le |\upsilon_d|^{d-s}
\le |\upsilon_d|^{d-1}
$$
and 
$$
 \rmM(\upsilon^{-1})=|\upsilon_1 \cdots \upsilon_s|^{-1}=\rmM(\upsilon)=\lambda
$$ 
with
$$ 
\lambda \le |\upsilon_1|^{-s}
\le |\upsilon_1|^{-(d-1)}.
$$ 
Therefore we have
\begin{equation}\label{Equation:hauteurunite}
\lambda^{1/(d-1)}\le |\upsilon_d|\le \lambda \quad\hbox{and}\quad
\lambda^{-1}\le |\upsilon_1| \le \lambda^{-1/(d-1)}. 
\end{equation}

\section{Examples}\label{S:example}
The lower bound $ \mu \ge \lambda^{\cst{kappa:maja3}}$ quoted in section $\ref{S:introduction}$ is true 
\\
$\bullet$ when $d=3$ and the cubic field $K$ is not totally real;
\\
$\bullet$ for the simplest fields of degree $3$ (see \cite{LW-Thue}), and also for the simplest fields of degrees $4$ and $6$;
\\
$\bullet$ when $-\upsilon$ is a Galois conjugate of $\upsilon$ (which means that the irreducible polynomial of $\upsilon$ is in $\Z[X^2]$), and more generally when $|\upsilon_1|=|\upsilon_2|$ and $|\upsilon_{d-1}|=|\upsilon_d|$ with $d\ge 4$.

Here is an example of this last situation. 
Let $\epsilon$ be an algebraic unit, not a root of unity, of degree $\ell\ge 2$ and conjugates $\epsilon_1,\epsilon_2,\dots,\epsilon_\ell$. Let $h\ge 2$ and let $d=\ell h$. For $a\in\Z$, define 
\begin{equation}\label{Equation:corollary}
F_a(X,Y)= \prod_{i=1}^\ell (X^h-\epsilon_i^a Y^h).
\end{equation}
Let $R$ be the regulator of the field $\Q(\epsilon^{1/h})$. 

From Theorem $\ref{Theorem:majorationdea}$ we deduce the following corollary. 

\begin{corollary}\label{Corollary:exemple}
Let $m\ge 1$. 
If the form $F_a$ in {\rm ($\ref{Equation:corollary}$)} is irreducible and if there exists $(x,y)\in\Z^2$ with $xy\not=0$ and $|F_a(x,y)|\le m$, then 
$$
 |a| 
\le \Newcst{kappa:exemple} (R+\log m + \log \house{\epsilon}) R \log^\star(R\log \house{\epsilon}).
$$
\end{corollary}

  {\sc Proof}.
Without loss of generality, assume $|\epsilon_1|\le |\epsilon_2|\le \dots\le |\epsilon_\ell|$, so that $|\epsilon_\ell|=\house{\epsilon}
$. Let $\zeta$ be a primitive $h$-th root of unity. Let $\upsilon=\epsilon^{1/h}$.
We apply Theorem $\ref{Theorem:majorationdea}$ with $\alpha=\zeta$, $a_0=1$, $\lambda_0=1$, 
$\lambda\le \house{\epsilon}^\ell$, $F_0(X,Y)=(X^h-Y^h)^\ell$ and
$$
\upsilon_{ih+j}=\zeta^{j-1} \epsilon_{i+1}^{1/h} \quad (0\le i\le \ell-1,\; 1\le j\le h).
$$
From $|\upsilon_1|=|\upsilon_2|=|\epsilon_1|^{1/h}<1$ and $|\upsilon_{d-1}|=|\upsilon_d|=|\epsilon_\ell|^{1/h}
$ we deduce 
$$
 \mu=
 \left|\frac{\epsilon_\ell}{\epsilon_1}\right|^{1/h}=\left|\frac{\upsilon_d}{\upsilon_1}\right| 
 $$
and using ($\ref{Equation:hauteurunite}$) we conclude
 $$
 \log\mu\ge\frac{2}{d-1} \log \lambda.   
 $$   \vskip -.9cm \hfill $\Box$
\mbox{}\smallskip

A variant of this proof is to take $\alpha=1$, $\lambda_0=1$, $F_0(X,Y)=(X-Y)^d$, and to use the fact that 
$\zeta^a$ is also a primitive $h$-th root of unity since $F_a$ is irreducible.

\section{Auxiliary results}

 \subsection{An elementary result}
 
 For the convenience of the reader, we include the following elementary result -- similar arguments are often used without explicit mention in the literature.
 
 \begin{lemma}\label{Lemma:elementary}
 Let $U$ and $V$ be positive numbers satisfying $\;U\le V\log^\star U$. Then \mbox{$U< 2V\log^\star V$.}
 \end{lemma}
 
 {\bf Proof.} 
 If $\log U \le 1$, the assumption is $U\le V$ and the conclusion follows. Assume $\log U>1$. Then 
 $ \log U\le \sqrt{U}$, hence the hypothesis of the lemma implies $U\le V\sqrt{U}$ and therefore we have $ U\le V^2$.
 We deduce 
 $$
 \log U\le 2\log V,
 $$
 hence
$$
U\le V\log U \le 2V\log V.
 $$
\vskip  -.8cm \hfill $\Box$
 
 \subsection{Diophantine tool}
In this section only, the positive integer $d$ is not restricted to $d\ge 3$. 

The main tool is the following Diophantine estimate (\cite[Proposition 2]{LW2}, \cite[Theorem 9.1]{GL326} or \cite[Th. 3.2.4]{EG}), the proof of which uses transcendental number theory. 

\begin{proposition}\label{Proposition:FormeLineaireLogarithmes}
Let $s$ and $D$ be two positive integers. There exists an effectively computable positive constant $\kappa(s,D)$, depending only upon $s$ and $D$, with the following property.
Let $\eta_1,\ldots, \eta_s$ be nonzero algebraic numbers generating a number field of degree $\le D$. Let $c_1,\ldots,c_s$ be rational integers and let $H_1,\ldots,H_s$ be real numbers $\ge 1$ satisfying 
$$
H_i\ge \rmh(\eta_i) \quad (1\le i\le s).
$$ 
 Let $C$ be a real number with $C\ge 2$. 
 Suppose that one of the following two statements is true:
 \\
\indent {\rm (i) } $ C\ge\max_{1\le j\le s} |c_j| $
\\
 or
 \\
\indent {\rm (ii)} 
 $H_j\le H_s$ for $1\le j\le s$ and 
 $$
 C\ge
\max_{1\le j\le s} \left\{
 \frac{H_j}{H_s} |c_j|\right\}. 
$$ 
Suppose also $\eta_1^{c_1}\cdots \eta_s^{c_s} \not=1$.
Then 
$$
|\eta_1^{c_1}\cdots \eta_s^{c_s}-1|>
\exp\{-
\kappa(s,D) H_1\cdots H_s\log C\}. 
$$
 \end{proposition}
 
The statement (ii) of Proposition $\ref{Proposition:FormeLineaireLogarithmes}$ implies the statement (i) by permuting the indices so that $H_j\le H_s$ for $1\le j\le s$; however, we find it more convenient to use the part (i) so that we can use the estimate without permuting the indices. 
 
 We will use Proposition $\ref{Proposition:FormeLineaireLogarithmes}$ several times. Here is a first consequence. 
 
 \begin{corollary}\label{Corollaire:DiophantineTool}
 Let $d\ge 1$. There exists a constant $\Newcst{kappa:corollairediophantien}$, which depends only on $d$, with the following property. Let $K$ be a number field of degree $d$. 
 Let $\alpha_1$, $\alpha_2$, $\upsilon_1$, $\upsilon_2$ be nonzero elements in $K$ and let $a$ be a nonzero integer. Set $\gamma_1=\alpha_1\upsilon_1^a$ and $\gamma_2=\alpha_2\upsilon_2^a$. Let $\lambda_0 $ and $\lambda$ satisfy 
 $$
 \max\{
 \rmh(\alpha_1), \rmh(\alpha_2)\}\le \log \lambda_0,
 \quad
 \max\{
 \rmh(\upsilon_1), \rmh(\upsilon_2)\}\le \log \lambda
 $$
 and assume $\gamma_1\not= \gamma_2$. Define
$$
\chi=
 (\log^\star \lambda_0)(\log^\star \lambda)\log^\star\left(|a|
 \min\left\{1,\; 
 \frac{\log^\star \lambda}{\log^\star \lambda_0}
 \right\}
 \right).
 $$ 
 Then 
$$
| \gamma_1- \gamma_2|\ge \max
\left\{
|\gamma_1|, |\gamma_2|\right\} e^{
-\cst{kappa:corollairediophantien}\chi}.
$$
 \end{corollary}
 
 {\bf Proof.} 
  By symmetry, \,without loss of generality, \, we may assume \mbox{$|\gamma_2|\ge |\gamma_1|$.} 
 Set
 $$
 s=2,\quad
 \eta_1=\frac{\upsilon_1}{\upsilon_2},
 \quad \eta_2=\frac{\alpha_1}{\alpha_2},
 \quad
 c_1=a,\quad c_2=1,
 $$
 $$
 H_1=2\log^\star\lambda,
 \quad 
 H_2= 2\log^\star\lambda_0,
 \quad
 C=\max \left\{2, |a| \min\left\{1, \frac{ H_1}{H_2}\right\}\right\} \cdot 
 $$
 The conclusion of Corollary \ref{Corollaire:DiophantineTool} follows from Proposition $\ref{Proposition:FormeLineaireLogarithmes}$ (via part (i) if $H_1\ge H_2$, via part (ii) otherwise), thanks to the relation
 $$
 \left|\eta_1^{c_1}\eta_2^{c_2}-1\right|=|\gamma_2|^{-1} | \gamma_1- \gamma_2|.
$$ 
 \vskip -.8cm \hfill $\Box$
 
 \subsection{Lower bound for the height and the regulator}

For the record, we quote Kronecker's Theorem and its effective improvement.

\begin{lemma}\label{Lemma:Kronecker} \,
{\rm (a)}
If a nonzero algebraic integer $\alpha$ has all its conjugates in the closed unit disc $\{ z\in \C \; \mid \; |z| \leq 1 \}$, then $\alpha$ is a root of unity. \\
\indent {\rm (b)}
More precisely, given $d\ge 1$, there exists an effectively computable positive constant $\Newcst{kappa:Dobrowolski}$, depending only on $d$, such that, if $\alpha$ is a nonzero algebraic integer of degree $d$ satisfying $\rmh(\alpha)<\cst{kappa:Dobrowolski}$, then $\alpha$ is a root of unity. 
\end{lemma}

{\bf Proof.}
Voutier (1996) refined an earlier estimate due to Dobrowolski (1979) by proving that the conclusion of the part (b) in Lemma $\ref{Lemma:Kronecker}$ holds with 
$$
\cst{kappa:Dobrowolski}=\begin{cases}
\log 2 &\mbox{ if } d=1,
\\
\displaystyle {2}{d(\log d)^3} &\mbox{ if }  d\ge 2.
\end{cases}
$$
See for instance \cite[Prop. 3.2.9]{EG} and \cite[\S 3.6]{GL326}. 
  \hfill $\Box$ 

\begin{lemma}\label{Lemma:Regulator}
There exists an explicit absolute constant $\Newcst{kappa:Regulator}>0$ such that the regulator $R$ of any number field of degree $\ge 2$ satisfies $R> \cst{kappa:Regulator}$. 
\end{lemma}

{\bf Proof.}
According to a result of Friedman (1989 -- see \cite[(1.5.3)]{EG}) the conclusion of Lemma $\ref{Lemma:Regulator}$ holds with $\cst{kappa:Regulator}= 0.2052$.
 \hfill $\Box$

\subsection{A basis of units of an algebraic number field}

Here is Lemma 1 of 
\cite{BG2debut}. See also \cite[Proposition 4.3.9]{EG}. The result is essentially due to C.L.~Siegel \cite{MR0249395}. 

\begin{proposition}\label{Proposition:hauteurbasegroupeunites}
Let $d$ be a positive integer with $d\ge 3$. There exist effectively computable constants $\Newcst{kappa:c_7},\Newcst{kappa:c_8},\Newcst{kappa:c_9}$ depending only on $d$, with the following property. 
Let $K$ be a number field of degree $d$, with unit group of rank $r$. Let $R$ be the regulator of this field. Denote by $\varphi_1,\varphi_2,\dots,\varphi_r$ a set of $r$ embeddings of $K$ into $\C$ containing the real embeddings and no pair of conjugate embeddings. Then there exists a fundamental system of units $\{\epsilon_{1},\epsilon_{2},\dots,\epsilon_{r}\}$ of $K$ which satisfies the following:

  {\rm (i)} $\displaystyle
\prod_{1\le i\le r} \rmh(\epsilon_{i})\le \cst{kappa:c_7} R$;

 {\rm (ii)} $\displaystyle \max_{1\le i\le r} \rmh(\epsilon_i)\le \cst{kappa:c_8} R$;

  {\rm (iii)} The absolute values of the entries of the inverse matrix of 
$$
(\log|\varphi_j(\epsilon_i)|)_{1\le i, j \le r}
$$ 
do not exceed $ \cst{kappa:c_9}$. 
\end{proposition}

 
 The next result is \cite[Lemma A.15]{ST}. 
 
 \begin{lemma}\label{Lemma:A15}
 Let $\epsilon_1$, $\epsilon_2$, $\dots$, $\epsilon_r$ be an independent system of units for $K$ satisfying the condition $(ii)$ of Proposition $\ref{Proposition:hauteurbasegroupeunites}$. Let $\beta\in\Z_K$ with $\rmN_{K/\Q}(\beta)=m\not=0$. Then there exist $b_1,b_2,\dots,b_r$ in $\Z$ and $\betatilde\in\Z_K$ with conjugates $\betatilde_1,\betatilde_2,\dots,\betatilde_d$, satisfying 
 $$
 \beta=
\betatilde \epsilon_1^{b_1}\epsilon_2^{b_2}\cdots\epsilon_r^{b_r}
 $$
 and
 $$ 
|m|^{1/d} e^{-\Newcst{kappa:A15} R}
 \le |\betatilde_j|
 \le |m|^{1/d} e^{\cst{kappa:A15} R}
 \quad \hbox{for} \quad j=1,\dots,d.
 $$
 \end{lemma}
 
 The conclusion of Lemma $\ref{Lemma:A15}$ can be written
 $$ 
\left|
\log
\left(
|m|^{-1/d} |\betatilde_j|
\right)
\right|
\le 
\cst{kappa:A15} R
 \quad \hbox{for} \quad j=1,\dots,d.
 $$
 
\subsection{Estimates for the conjugates}\label{SS:Conjugates} 

\begin{lemma}\label{Lemma:deuxconjuguesreels}
Let $\gamma$ be an algebraic number of degree $d\ge 3$. Let 
$\gamma_1$,\linebreak 
$\gamma_2$, $\dots$, $\gamma_d$  be the conjugates of $\gamma$ with $ |\gamma_1|\le |\gamma_2|\le \cdots \le |\gamma_d|$. 
\\
\indent {\rm (a)} If $|\gamma_1|< |\gamma_2|$ and $\gamma_2\in\R$, then $|\gamma_2|<|\gamma_3|$. 
\\
\indent {\rm (b)} If $|\gamma_{d-1}|< |\gamma_d|$ and $\gamma_{d-1}\in\R$, then $|\gamma_{d-2}|<|\gamma_{d-1}|$. 
\end{lemma} 
{\bf Proof.}
(a) The conditions $|\gamma_1|< |\gamma_2|\le |\gamma_i|$ for $3\le i\le d$ imply that $\gamma_1$ is real and that $-\gamma_1$ is not a conjugate of $\gamma_1$. Hence the minimal polynomial of $\gamma$ is not a polynomial in $X^2$. Assume $|\gamma_2|=|\gamma_3|$. Since $-\gamma_2$ is not a conjugate of $\gamma_2$, we deduce $\gamma_3\not\in\R$, hence $d\ge 4$. We may assume $\gamma_4=\overline{\gamma_3}$. Let $\sigma$ be an automorphism of $\Qbar$ which maps $\gamma_2$ to $\gamma_1$; via $\sigma$, let $\gamma_j$ be the image of $\gamma_3$ and $\gamma_k$ the image of $\gamma_4$. From 
$$
\gamma_2^2=\gamma_3\gamma_4
$$
we deduce $
\gamma_1^2=\gamma_j\gamma_k
$
and
$
|\gamma_1|^2=|\gamma_j\gamma_k|$. 
This is not possible since $|\gamma_j|>|\gamma_1|$ and $|\gamma_k|>|\gamma_1|$. 

\indent (b) We deduce (b) from (a), by using $\gamma\mapsto 1/\gamma$ (or by repeating the proof, {\it mutatis mutandis}). 
\hfill $\Box$

{\bf Remark.} Here is an example showing that the assumptions of Lemma $\ref{Lemma:deuxconjuguesreels}$ are sharp. 
The polynomial $X^4-4X^2+1$ is irreducible, its roots are 
$$
\upsilon_1=\sqrt{2-\sqrt{3}},\quad \upsilon_2=-\upsilon_1,\quad \upsilon_3=1/\upsilon_1=\sqrt{2+\sqrt{3}},\quad \upsilon_4=-\upsilon_3
$$
with 
$$
\upsilon_1=|\upsilon_2|<\upsilon_3=|\upsilon_4|.
$$
More generally, if $h\geq 2$ is a positive integer and $\epsilon$ is a quadratic unit with Galois conjugate $\epsilon'$ and if $\epsilon^{1/h}$ has degree $2h$, then it has $h$ conjugates of absolute value $|\epsilon|^{1/h}$ and $h$ conjugates of absolute value $|\epsilon'|^{1/h}$.
See also \S$\ref{S:example}$. 

\begin{lemma}\label{Lemma:minorationalpha2suralphad}
Let $\upsilon$ be an algebraic unit of degree $d\ge 3$. Set $\lambda=\rmM(\upsilon)$. Let $\upsilon'$ and $\upsilon''$ be two conjugates of $\upsilon$ with $ |\upsilon'| <|\upsilon''|$. Then 
$$
\log\frac{|\upsilon''|}{ |\upsilon'|}\ge \Newcst{kappa:lemme1} \lambda^{-(d^3+2d^2-d+2)/2}.
$$
\end{lemma}

We will deduce Lemma $\ref{Lemma:minorationalpha2suralphad}$ from Theorem 1 of \cite{GS} which\footnote{This reference was kindly suggested to us by Yann Bugeaud.} states the following. 

\begin{lemma}[X.~Gourdon and B.~Salvy \cite{GS}]\label{Lemma:GourdonSalvy}
Let $P$ be a polynomial of degree $d\ge 2$ with integer coefficients and with Mahler measure $\rmM(P)$. If $\alpha'$ and $\alpha''$ are two roots of $P$ with $ |\alpha'| <|\alpha''|$, then 
$$
|\alpha''| - |\alpha'|\ge \Newcst{kappa:lemmeGourdonSalvy} M(P)^{-d(d^2+2d-1)/2}
$$
with 
$$
\cst{kappa:lemmeGourdonSalvy}=\frac{\sqrt{3}}{2} \bigl( d(d+1)/2\bigr) ^{-d(d+1)/4-1}.
$$
\end{lemma}

{ \bf Proof of Lemma $\ref{Lemma:minorationalpha2suralphad}$}\,.
We apply Lemma $\ref{Lemma:GourdonSalvy}$ to the minimal polynomial of $\upsilon$. 
To conclude the proof of Lemma $\ref{Lemma:minorationalpha2suralphad}$, we use the bounds $|\upsilon'|\le \lambda$ and 
$$
 \log (1+x)\ge \frac{x}{2} \quad\hbox{for}\quad 0\le x\le 1
\quad\hbox{with} \quad
x= \frac{|\upsilon''|}{ |\upsilon'|} -1.
$$
\vskip -1cm \hfill $\Box$

\section{Proof of Theorem {\bf \ref{Theorem:majorationdea}} } \label{S:ProofTheorem2}

Theorem $\ref{Theorem:majorationdea}$ with the assumption $|F_a(x,y)| \leq m$ will be secured if we deal with the equation $F_a(x,y) = m$ 
with $m\neq 0$.

Let $(a,x,y,m)\in \Z^4$ satisfy $m\not=0$, $xy\not=0$, $[\Q(\alpha\upsilon^a):\Q]=d$ and 
$$
F_a(x,y)=m.
$$ 
Without loss of generality, we may restrict $(a,y)$ to $a\ge 0$ (otherwise, replace $\upsilon$ by $\upsilon^{-1}$) and to $y>0$ (otherwise replace $F_a(X,Y)$ by $F_a(X,-Y)$). 
 
The form $\tilde{F}_a(X,Y)=a_0^{d-1}F_a(X,Y)$ has coefficients in $\Z$, and if we set $\tilde{x}=a_0x$, $\tilde{y}=y$, $\tilde{m}=a_0^{d-1}m$ we have 
$\tilde{F}_a(\tilde{x},\tilde{y})=\tilde{m}$ with $(\tilde{x},\tilde{y})\in\Z^2$. Therefore, there is no loss of generality to assume $a_0=1$.

Theorem $\ref{Theorem:majorationdea}$ includes the assumption that $\upsilon$ is not a root of unity, hence $\lambda>1$. More precisely, it follows from the part (b) of Lemma $\ref{Lemma:Kronecker}$ that 
$$
\log \lambda\ge \cst{kappa:Dobrowolski}.
$$
In particular, we have 
 $$
 \log^\star \lambda\le \max\left\{
 1, \frac{1}{\cst{kappa:Dobrowolski}}
 \right\}
 \log \lambda,
 $$
an inequality which can be written 
\begin{equation}\label{Equation:logstarlambda}
 \log^\star \lambda\le \cst{lambda}\log \lambda
 \end{equation}
 with an effectively computable constant $ \Newcst{lambda}>0$. 
 
From Lemma $\ref{Lemma:Regulator}$, we deduce that $R>\cst{kappa:Regulator}$.
Therefore, there is no loss of generality to assume that, for a sufficiently large constant $\Newcst{kappa:minorationdea}$, we have
\begin{equation}\label{Equation:minorationdea}
a \ge \cst{kappa:minorationdea} \bigl(\log |m| + (\log^\star \lambda_0)\log^\star\log^\star\lambda\bigr).
 \end{equation}
 This hypothesis will frequently be used, sometimes without explicit mention.
 
 By assumption, $\Q(\alpha \upsilon^a)=K$. If some conjugate $\sigma_j(\alpha \upsilon^a)$ of $\alpha \upsilon^a$ is real, then it follows that $\sigma_j(K)\subset\R$, hence the embedding $\sigma_j$ is real, and $\alpha_j$ and $\upsilon_j$ are both real. We also notice that if $\sigma_j(\upsilon)=-\sigma_i(\upsilon)$ with $i\not=j$,
then it follows that $\upsilon$ and $-\upsilon$ are conjugate, hence the irreducible polynomial of $\upsilon$ belongs to $\Z[X^2]$.

Recall that $\upsilon_i=\sigma_i(\upsilon)$ ($i=1,\dots,d$) and that
$$
 |\upsilon_1|\le |\upsilon_2|\le \cdots \le |\upsilon_d|. 
 $$
Let us write $\alpha_i$ for $\sigma_i(\alpha)$ ($i=1,\dots,d$).
 Let 
 $$
 \gamma=\alpha \upsilon^a
 \quad
 \hbox{and}\quad
\beta=x-\gamma y.
$$
Since $a_0=1$, it follows that $\alpha$, $\beta$ and $\gamma$ are algebraic integers in $K$. 
For $j=1,2,\dots,d$, define $\gamma_j$ and $\beta_j$ by
$$
\gamma_j=\sigma_j(\gamma)=\alpha_j \upsilon_j^a,\quad
\beta_j=\sigma_j(\beta)=x-\alpha_j\upsilon_j^a y=x-\gamma_j y.
$$
The assumption $F_{a}(x,y)=m$ yields $\beta_1\beta_2\cdots \beta_d=m$. Let $i_0\in \{1,2,\dots,d\}$ be an index such that 
$$
|\beta_{i_0}|=\min_{1\le i\le d} |\beta_i|.
$$

We define $\Psi_1,\Psi_2,\dots,\Psi_d$ by the following conditions:
 $$
 \beta_i=
\begin{cases}
\gamma_{i_0}y \Psi_i
&\hbox{for $1\le i<i_0$},
 \\
\gamma_iy \Psi_i
&\hbox{for $i_0< i\le d$}
 \end{cases}
$$
and
$$
 \beta_{i_0} = \frac{m}{y^{d-1} } \cdotp \frac{\gamma_1\gamma_2\cdots \gamma_{i_0-1}}{
\gamma_{i_0}^{i_0-2}}
\Psi_{i_0}  \cdot
 $$
 We split the proof into several steps. 

\indent
{\tt Step 1.}   We start by proving that
 \begin{equation}\label{Eq:majx}
 |x|\le 2\lambda_0\lambda^a y
\end{equation}
and that there exists an effectively computable positive constant $\Newcst{kappa:step1}$ depending only on $d$ such that
 \begin{equation}\label{Eq:majPsi}
 e^{-\cst{kappa:step1} \chi} \le |\Psi_i| \le e^{\cst{kappa:step1} \chi} \qquad (i=1,2,\dots,d)
\end{equation}
with 
$$
\chi=
 (\log^\star \lambda_0)(\log
 \lambda)\log
 \left( a
 \min\left\{1,\;
 \frac{\log
 \lambda}{\log^\star \lambda_0}
 \right\}
 \right).
 $$ 
From the estimate $(\ref{Eq:majPsi})$ we will deduce
$$
| \beta_{i_0} |<|\beta_i|
$$ 
for $i\not=i_0$, which implies $\alpha_{i_0}\in\R$ and $\upsilon_{i_0}\in \R$. 

\indent
{\tt Remark.} The estimate $(\ref{Eq:majPsi})$ can be written as follows: 
$$
\left|
 \log \left(
 |\beta_i|y^{-1}|
\max\{|\gamma_i^{-1}|, |\gamma_{i_0}^{-1} |\}
\right)
 \right|
 \le \cst{kappa:step1}\chi
$$
for $i\not=i_0$ and 
$$
 \left|
 \log \left(
 |\beta_{i_0}|
 \frac{y^{d-1} }{|m|}
\left | \gamma_1^{-1} \cdots \gamma_{i_0-1}^{-1} \gamma_{i_0}^{i_0-2}
 \right|
 \right)
 \right|
 \le \cst{kappa:step1}\chi.
 $$
 
{\tt Proof}  of $(\ref{Eq:majx})$ and $(\ref{Eq:majPsi})$.
We have
\begin{equation}\label{Eq:premieremajorationx}
|x|=|\beta_{i_0}+\gamma_{i_0}y| \le |\beta_{i_0}|+|\gamma_{i_0}|y.
\end{equation}
From $|\beta_{i_0}|\le |\beta_i|$ for $i=1,2,\dots,d$ and $\beta_1\cdots\beta_d=m$, we deduce \linebreak
$|\beta_{i_0}|\le |m|^{1/d}$, hence
$$
|x| \le |m|^{1/d}+|\gamma_{i_0}|y\le |m|^{1/d}+\lambda_0\lambda^a y.
$$
Using the assumption $(\ref{Equation:minorationdea})$, we check $|m|^{1/d}\le \lambda_0\lambda^a y$, whereupon the inequality $(\ref{Eq:majx})$ is secured.
\par

We also have
\begin{equation}\label{Eq:majorationbetai0}
|\beta_{i_0}|^{d-1}\max_{1\le i\le d} |\beta_i|\le |m|.
\end{equation}
For $i=1,2,\dots,d$, we write 
\begin{equation}\label{Eq:betai}
\beta_i=\beta_{i_0}+y(\gamma_{i_0}-\gamma_i).
\end{equation}

We have
$$
|\alpha_1\alpha_2\cdots\alpha_d|\ge 1
$$
(recall $a_0=1$), 
hence
$$
\frac{1}{\lambda_0 }\le |\alpha_i|\le \lambda_0
\quad\hbox{
for $i=1,2,\dots,d$. }
$$

We choose an index $j_0\not=i_0$ as follows: 

$\bullet$
If $|\upsilon_{i_0}|\le \lambda^{1/(2(d-1))}$, we take $j_0=d$ so that, with the help of $(\ref{Equation:hauteurunite}$), we have $|\upsilon_{j_0}|\ge \lambda^{1/(d-1)}$, whereupon with the help of $(\ref{Equation:minorationdea})$ we obtain
$$
\left|
\frac{\gamma_{i_0}}{\gamma_{j_0}}
\right|<\frac{1}{2}\cdotp
$$

$\bullet$
If $|\upsilon_{i_0}|> \lambda^{1/(2(d-1))}$, we take $j_0=1$ so that, again with the help of $(\ref{Equation:hauteurunite}$), we have $|\upsilon_{j_0}|\le \lambda^{-1/(d-1)}$, whereupon with the help of $(\ref{Equation:minorationdea})$ we obtain
$$
\left|
\frac{\gamma_{j_0}}{\gamma_{i_0}}
\right|<\frac{1}{2}\cdotp
$$

In both cases, we deduce 
$$
|\gamma_{j_0}-\gamma_{i_0} | \ge
\frac{1}{2}\max
\{
|\gamma_{j_0}|, |\gamma_{i_0} | 
\}
\ge \frac{\lambda^{a/(2(d-1))}}{2\lambda_0 }
$$
and therefore, using $(\ref{Equation:minorationdea})$ again together with $(\ref{Eq:majorationbetai0})$ and $(\ref{Eq:betai})$, we obtain 
$$
|\beta_{j_0}|\ge 
|\gamma_{j_0} -\gamma_{i_0} | y -|\beta_{i_0}| \ge
\frac{\lambda^{a/(2(d-1))}y} {2\lambda_0} -|m|^{1/d} \ge
\lambda^{a/(2d)} y.
$$
Since $\displaystyle\max_{1\leq i\leq d} |\beta_i| \geq \lambda^{a/(2d)} y$, from 
$(\ref{Eq:majorationbetai0})$ we deduce 
\begin{equation}\label{Eq:step1a}
|\beta_{i_0}|\le \left(\frac{|m|}{y \lambda^{a/(2d)}}\right)^ {1/(d-1)}.
\end{equation} 
In particular, thanks to $(\ref{Equation:minorationdea})$, we have 
\begin{equation}\label{Eq:betoi0le12}
|\beta_{i_0}|\le \frac{1}{2}\cdotp
\end{equation}
Using the assumption $|x|\ge 1$ together with $(\ref{Eq:premieremajorationx})$, we deduce 
\begin{equation}\label{Eq:minorationalphai0}
\frac{|x|}{2} \leq |\gamma_{i_0}| y \leq |x]+|\beta_{i_0}| \leq \frac{3|x|}{2} 
\cdotp
\end{equation}
\par

Let $i\not=i_0$. The upper bound 
$$
|\gamma_i-\gamma_{i_0}|\le 2
\max\{ |\gamma_{i_0}|, |\gamma_i|\}
$$
is trivial, while the lower bound 
\begin{equation}\label{equation:minorationalphaimoinsalphai0}
|\gamma_i-\gamma_{i_0}|\ge 
\max\{ |\gamma_{i_0}|, |\gamma_i|\}
e^{-\cst{kappa:minorationalphaimoinsalphai0}\chi}
\end{equation} 
follows from $(\ref{Equation:logstarlambda})$ and from Corollary \ref{Corollaire:DiophantineTool}.
We first use the lower bound 
$$
|\gamma_i-\gamma_{i_0}|
\ge 
|\gamma_{i_0}|
e^{-\Newcst{kappa:minorationalphaimoinsalphai0}\chi}. 
$$
Using $(\ref{Eq:minorationalphai0})$, we obtain 
\begin{equation}\label{Eq:alphaia}
|\gamma_i-
\gamma_{i_0}|
\ge
\frac{1}{2y}
e^{-\cst{kappa:minorationalphaimoinsalphai0}\chi}\ge
\frac{2}{y}
e^{-\cst{kappa:minorationalphaimoinsalphai0bis}\chi}
\end{equation} 
with $\Newcst{kappa:minorationalphaimoinsalphai0bis}>0$.
Using the contrapositive of Lemma $\ref{Lemma:elementary}$ with 
$$
U=a\frac{\log^\star \lambda}{\log^\star \lambda_0}, \quad 
V=\frac{1}{\Newcst{kappa:majorationchi} }\log^\star \lambda,
$$
we deduce from $(\ref{Equation:minorationdea})$ that
$$
\chi\le \cst{kappa:majorationchi} a \log ^\star \lambda.
$$
Recall that $\cst{kappa:minorationdea}$ is sufficiently large, hence $\cst{kappa:majorationchi}$ is sufficiently small. 
Now from $(\ref{Eq:step1a})$,
the inequality $|m|\leq e^{|a]/\cst{kappa:minorationdea}}$ 
and $(\ref{Eq:alphaia})$ we deduce
$$
|\beta_{i_0}|\le 
 |m| ^{1/(d-1)} \lambda^{-a/(2d(d-1)) }
 \le
 \lambda^{-\Newcst{kappa:minorationalphaimoinsalphai04}a}\le 
e^{-\cst{kappa:minorationalphaimoinsalphai0bis}\chi}\le 
\frac{1}{2} y |\gamma_i-\gamma_{i_0}|.
$$
Therefore, for $i\not=i_0$, 
using $(\ref{Eq:betai})$, we deduce 
$$ 
\frac{1}{2} y |\gamma_i-\gamma_{i_0}|
\le |\beta_i|\le \frac{3}{2} y |\gamma_i-\gamma_{i_0}|.
$$ 
Using once more $(\ref{equation:minorationalphaimoinsalphai0})$, we obtain $(\ref{Eq:majPsi})$ for $i\not=i_0$. We also deduce 
\begin{equation}\label{Equation:minorationbetai}
|\beta_i|>\lambda^{-a/(2d)}
\quad \hbox{for}\quad i\not=i_0.
\end{equation}

Recall 
$$
\rmN(\gamma)=
\gamma_1\gamma_2\cdots\gamma_d=
\rmN(\alpha)\rmN(\upsilon)^a=\pm \rmN(\alpha)
\quad\hbox{and}\quad 
\rmN(\beta)=\beta_1\beta_2\cdots\beta_d=m.
$$
The estimate $(\ref{Eq:majPsi})$ for $i=i_0$ follows from the relations
$$
\Psi_1\Psi_2\cdots\Psi_d\rmN(\gamma)=1,
$$
$$
\frac{m}{\beta_{i_0}}=\prod_{i\not=i_0} \beta_i
=
y^{d-1}\gamma_{i_0}^{i_0-1}\gamma_{i_0+1} \cdots \gamma_d\prod_{i\not=i_0} \Psi_i
$$
and
$$
\frac{\rmN(\gamma)}{
\gamma_{i_0}^{i_0-1}
\gamma_{i_0+1}\cdots\gamma_d}=
\frac{\gamma_1\cdots\gamma_{i_0-1}}{
\gamma_{i_0}^{i_0-2}} \cdotp
$$ 
From $(\ref{Equation:minorationdea})$ and $(\ref{Eq:majPsi})$, 
we deduce 
$$
|\beta_{i_0}|< \frac{|m|}{y^{d-1}}
|\gamma_1| e^{\cst{kappa:step1}\chi}
<\lambda^{-a/(2d)},
$$
hence from $(\ref{Equation:minorationbetai})$ we infer 
$|\beta_{i_0}|< |\beta_i|$ for $i\not=i_0$. It follows that $\beta_{i_0}$ is real, and therefore $\gamma_{i_0}$, $\alpha_{i_0}$ and $\upsilon_{i_0}$ also. 
\hfill $\Box$

\indent
{\tt Step 2}. Let $\epsilon_1,\epsilon_2,\dots,\epsilon_r$ be a basis of the group of units of $K$ given by Proposition $\ref{Proposition:hauteurbasegroupeunites}$. From Lemma $\ref{Lemma:A15}$, it follows that there exists $\betatilde\in\Z_K$ and $b_1,b_2,\dots,b_r$ in $\Z$ with
 $$
 \beta =
\betatilde\epsilon_1^{b_1}\epsilon_2^{b_2}\cdots\epsilon_r^{b_r}
 $$
 and
$$
|m|^{1/d} e^{-\cst{kappa:A15} R}\le
|\betatilde_i|
 \le |m|^{1/d} e^{\cst{kappa:A15} R}
 \quad\hbox{for}\quad i=1,2,\dots,d.
 $$
We set 
\begin{equation}\label{Equation:DefinitionB}
B= \cst{kappa:B} (R+a\log \lambda+\log y 
)
\end{equation}
with a sufficiently large constant $\Newcst{kappa:B} $.
We want to prove that
$$
\max_{1\le i\le r} |b_i|\le B.
$$ 

{\tt Proof.}
We consider the system of $d$ linear forms in $r$ variables with real coefficients
$$
L_i(X_1,X_2,\dots,X_r)= \sum_{j=1}^r X_j\log |\sigma_i(\epsilon_j)|,\quad (i=1,2,\dots,d).
$$
The rank is $r$. By Proposition $\ref{Proposition:hauteurbasegroupeunites}$(ii),
$$
\log |\sigma_i(\epsilon_j)|\le \Newcst{kappa:prop} R.
$$
For $i=1,2,\dots,d$, define $e_i=L_i(b_1,b_2,\dots,b_r)$. We have 
$$e_i=\log |\sigma_i(\beta/\betatilde)|=\log |\beta_i/\betatilde_i|,$$
 hence, using 
the inequality $|m|\leq e^{|a]/\cst{kappa:minorationdea}}$ and $(\ref{Eq:majPsi})$, 
we deduce
$$
|e_i|\le \Newcst{kappa:ei} (R+a\log \lambda+\log y 
).
$$
Computing $b_1,b_2,\dots,b_r$ by means of the system of linear equations
$$
L_i(b_1,b_2,\dots,b_r)=e_i \quad (i=1,2,\dots,d)
$$
and using Proposition $\ref{Proposition:hauteurbasegroupeunites}$(iii), 
we deduce
$$
\max_{1\le j\le r}|b_j|\le \Newcst{kappa:bj1} \max_{1\le i\le d} |e_i|
\le B. 
$$
\vskip -1cm \hfill $\Box$\bigskip

{\tt Step 3}. 
From the inequality (3.2) in \cite[Theorem 3]{BG2debut} (see also \cite[Th.~9.6.2]{EG}), thanks to $(\ref{Equation:minorationdea})$, we deduce the following upper bound for $|x|$ and $|y|$ in terms of $a$, $\lambda$, $\lambda_0$, $m$ and $R$: there exists a positive effectively computable constant $ \Newcst{kappa:majorationdexety}$ depending only on $d$ such that 
\begin{equation}\label{Equation:majorationdexety}
\log \max \{ |x|, y \}
\le \cst{kappa:majorationdexety} 
R (\log^\star R) \bigl(R+ 
a\log \lambda
\bigr).
\end{equation} 
 
\indent
{\tt Step 4}. 
Assume $c\gamma_i \beta_j \neq \gamma_k \beta_\ell$
for some indices $i,j,k,\ell$ in $\{1,\dots,d\}$ and some $c\in\{1,-1\}$.
Then there exists $ \Newcst{kappa:step5}>0$ such that 
\begin{align}\notag
\left|
c \frac{
\gamma_i\beta_j}{\gamma_k \beta_\ell}-1
\right|
\ge 
\exp\Biggl\{- \cst{kappa:step5} (\log \lambda)
&
(R+\log |m|+\log \lambda_0+\log \lambda)R
\\
\notag
&
\times \log 
\left( 
\frac{Ra\log\lambda}{R+\log |m|+\log\lambda_0+\log\lambda}
\right) 
\Biggr\}\cdotp
\end{align}

{\tt Proof.}
This lower bound follows from Proposition $\ref{Proposition:FormeLineaireLogarithmes}$(ii) with 
$$
\frac{
c\gamma_i\beta_j}{\gamma_k \beta_\ell}
=\eta_1^{c_1}\eta_2^{c_2}\cdots \eta_s^{c_s},
$$ where $s=r+2$ and
$$
\eta_t=\frac{\sigma_j(\epsilon_t)}{\sigma_\ell(\epsilon_t)}
\; (1\le t\le r),
\quad
\eta_{r+1}=\frac{\sigma_i(\upsilon)}{\sigma_k(\upsilon)},
\quad
\eta_{r+2}=\frac{c\sigma_j(\betatilde)\sigma_i(\alpha)}{\sigma_\ell(\betatilde)\sigma_k(\alpha)},
$$
$$
c_t=b_t \; (1\le t\le r), \quad c_{r+1}=a, \quad c_{r+2}=1,
$$ 
$$
H_t=\max\{1,2\rmh(\epsilon_t)\} \; (1\le t\le r),
$$
$$
H_{r+1}=\max\{1,2\log\lambda\},
\quad
 H_{r+2}=\Newcst{kappa:majbetatilde} (R+\log |m|+\log\lambda_0+\log\lambda),
$$
$$
C=2+\frac{2a\log\lambda +RB}{H_{r+2}}\cdotp
$$
Using Proposition $\ref{Proposition:hauteurbasegroupeunites} $(i) together with the part (b) of Lemma $\ref{Lemma:Kronecker}$, we deduce
$$
H_1H_2\cdots H_r\le\Newcst{kappa:Hi} R.
$$
Finally we deduce from the steps 2 and 3 that
$$
\log C\le \Newcst{kappa:logC}
 \log 
\left( 
\frac{Ra\log\lambda}{R+\log |m|+\log\lambda_0+\log\lambda}
\right),
$$
and this secures the above linear bound.
\hfill $\Box$

\indent
{\tt Step 5.} We will prove Theorem \ref{Theorem:majorationdea} by assuming 
$$
\max_{1\le i\le d}
\max\{|\Psi_i|, |\Psi_i|^{-1}\} > \mu^{ a/4}.
$$
Using $(\ref{Eq:majPsi})$, we deduce from our assumption
$$
\frac{ a}{4}\log \mu < \cst{kappa:step1} \chi,
$$
hence
$$
 a\le \frac{4\cst{kappa:step1}}{\log\mu}
 (\log^\star \lambda_0)(\log^\star \lambda)\log^\star\left( a
 \frac{\log^\star \lambda}{\log^\star \lambda_0}
 \right).
$$
With 
$$
 U=\frac{ a \log^\star \lambda}{\log^\star \lambda_0}
\quad \hbox{and} \quad
 V= \frac{4\cst{kappa:step1} (\log^\star \lambda)^2}{\log\mu} ,
 $$
 we have
 $U\le V\log^\star U$,
 and we conclude that 
we can use Lemma $\ref{Lemma:elementary}$ to deduce 
$$
 a\le
 \frac{8\cst{kappa:step1}(\log^\star \lambda_0)(\log^\star \lambda)}{\log\mu}
 \log \left(
 \frac{4\cst{kappa:step1}(\log^\star \lambda)^2}{\log\mu}
 \right),
$$
and the conclusion of Theorem \ref{Theorem:majorationdea}
 follows.
 
 In the rest of the paper, we assume
 \begin{equation}\label{equation:step5}
\max_{1\le i\le d}
\max\{|\Psi_i|, |\Psi_i|^{-1}\}\le \mu^{ a/4}.
\end{equation}
 
\indent
{\tt Step 6}. Our next goal is to prove the following results. 

\indent (a)
Assume $1\le i_0\le d-2$ and 
$$
\frac{|\upsilon_{d-1}|}{|\upsilon_{i_0}|}\ge \sqrt{ \mu}.
$$
Then 
$$ 
0<\left|
\frac{\gamma_{d-1}\beta_d}{\gamma_d \beta_{d-1}}-1
\right|
\le 4\lambda_0^2 \mu^{- a/4}.
$$
 
\indent (b) 
Assume $3\le i_0\le d$ and 
$$
\frac{|\upsilon_{i_0}|}{|\upsilon_2|}\ge \sqrt{ \mu}.
$$
Then
$$
0<
\left|
\frac{ \beta_1}{\beta_2}-1
\right| 
\le 2\lambda_0^2 \mu^{- a/4}.
$$
\indent (c) Assume $2\le i_0\le d-1$ and
$$
\min \left\{ \frac{ |\upsilon_{i_0}|}{ |\upsilon_1|}, \frac{ |\upsilon_d|}{ |\upsilon_{i_0}|}\right\}\ge \mu.
$$
Then
$$
\left|
\frac{
\gamma_{i_0}\beta_d}{\gamma_d \beta_1}+1
\right|
\le 4 |m| \lambda_0^4 \mu^{- a/2}.
$$

{\tt Proof} 
(a) We approximate $\beta_d$ by $-\gamma_dy$, $\beta_{d-1}$ by $-\gamma_{d-1}y$ and we eliminate $y$. Since $\gamma$ has degree $d$, we have
$$
\beta_d\gamma_{d-1}-\beta_{d-1}\gamma_d=x(\gamma_{d-1}-\gamma_d)\not=0.
$$
From 
$(\ref{Eq:minorationalphai0})$ we deduce $|x|\le 2|\gamma_{i_0} y|$ and
$$
|\beta_d\gamma_{d-1}-\beta_{d-1}\gamma_d|
\le
2 |\gamma_{i_0}| (|\gamma_d|+|\gamma_{d-1}|) y.
 $$
 Using $\beta_{d-1}=\gamma_{d-1}y\Psi_{d-1}$ together with the assumption
 $$
 |\upsilon_d|\ge|\upsilon_{d-1}|\ge \sqrt{ \mu} |\upsilon_{i_0}|,
 $$
 we deduce
 $$
 \left|
\frac{\gamma_{d-1}\beta_d}{\gamma_d \beta_{d-1}}-1
\right|
\le \frac{2 |\gamma_{i_0}|(|\gamma_{d-1}|+|\gamma_d|)}{|\gamma_{d-1}\gamma_d|}|\Psi_{d-1}|^{-1}\le
4\lambda_0^2\mu^{-a/2}|\Psi_{d-1}|^{-1}.
$$
The conclusion of (a) follows from ($\ref{equation:step5}$).

\indent
(b) We approximate $\beta_1$ and $\beta_2$ by $x$ and we eliminate $x$. Since $\gamma_1\not=\gamma_2$, we have 
$$
|\beta_1-\beta_2|= |(\gamma_2-\gamma_1) y| \not=0.
$$
From $\beta_2=\gamma_{i_0} y \Psi_2$ and the assumption
$$
 |\upsilon_1|\le|\upsilon_2|\le \mu^{-1/2} |\upsilon_{i_0}|,
 $$
we deduce 
$$
\left|
\frac{ \beta_1}{\beta_2}-1
\right|
\le
\frac{|\gamma_2|+|\gamma_1|}{|\gamma_{i_0}|} |\Psi_2|^{-1}
 \le 2\lambda_0^2\mu^{-a/2} |\Psi_2|^{-1}.
$$
Again, the conclusion of (b) follows from ($\ref{equation:step5}$).

\indent
(c) 
We approximate $\beta_1$ by $x$, $\beta_d$ by $-y\gamma_d$ and $x$ by $y \gamma_{i_0}$, then we eliminate $x$ and $y$. More precisely we have
$$
\beta_1\gamma_d+\beta_d\gamma_{i_0}=
(\gamma_d+\gamma_{i_0})\beta_{i_0} +\gamma_{i_0}^2y-\gamma_1\gamma_dy.
$$
 Hence 
$$
\frac{
\gamma_{i_0}\beta_d}{\gamma_d \beta_1}+1=
\frac{
(\gamma_d+\gamma_{i_0})\beta_{i_0}}{\gamma_d\beta_1} +
\frac{\gamma_{i_0}^2y}{\gamma_d\beta_1} -
\frac{\gamma_1 y}{ \beta_1}
\cdotp
$$ 
We have $\beta_1=\gamma_{i_0}y\Psi_1$. Therefore we have  
$$
\frac{|\gamma_{i_0}|^2y}{|\gamma_d\beta_1|}
=\frac{|\gamma_{i_0}|}{|\gamma_d|} |\Psi_1|^{-1} \le \lambda_0^2 \left|
\frac{\upsilon_{i_0}}{\upsilon_d}\right|^a |\Psi_1|^{-1}
$$
and                    
$$
\frac{|\gamma_1| y}{ |\beta_1|}=\frac{|\gamma_1| }{ |\gamma_{i_0}|} |\Psi_1|^{-1}
 \le
 \lambda_0^2 \left|
\frac{\upsilon_1}{\upsilon_{i_0}}\right|^a |\Psi_1|^{-1}.
$$
Finally, from 
$$
|\beta_{i_0}|\le \frac{|m|}{y^{d-1}}|\gamma_1\Psi_{i_0}|
$$ 
we deduce   
$$  \begin{array}{ll}
\displaystyle \frac{
(|\gamma_d|+|\gamma_{i_0}|)|\beta_{i_0}|}{|\gamma_d\beta_1|} 
&\le
(1+ \lambda_0^2)\left| \displaystyle 
\frac{\beta_{i_0}}{\beta_1}\right|\le
(1+ \lambda_0^2) \displaystyle 
\frac{|m|}{ y^d}\displaystyle 
\frac{|\gamma_1\Psi_{i_0}|}{|\gamma_{i_0}\Psi_1|}\\[6mm]
&\le
 \lambda_0^2(1+ \lambda_0^2)\displaystyle 
\frac{|m|}{ y^d}
 \left|
\frac{\upsilon_1}{\upsilon_{i_0}}
\right|^a\displaystyle 
\frac{|\Psi_{i_0}|}{ |\Psi_1|}\cdotp
\end{array}
$$
Hence from the assumptions  
$$ 
|\upsilon_1|\le
\mu^{-1} |\upsilon_{i_0}|
 \quad\hbox{and}
 \quad
 |\upsilon_{i_0}| \le
 \mu^{-1} |\upsilon_d|,
$$
we deduce  
$$
\left|
\frac{\gamma_{i_0}\beta_d}{\gamma_d \beta_1}+1
\right|
\le 4|m| \lambda_0^4 \mu^{- a } 
\frac{|\Psi_{i_0}|}{ |\Psi_1|}\cdot
$$
The conclusion of (c) follows from ($\ref{equation:step5}$).
\hfill $\Box$        

\indent
{\tt Step 7.} 
(a) Assume $|\upsilon_{i_0}|=|\upsilon_1|$. Since $\upsilon_{i_0}\in\R$, we deduce from Lemma $\ref{Lemma:deuxconjuguesreels}$ that
$$
|\upsilon_1|<|\upsilon_{d-1}|.
$$

If $|\upsilon_2|<|\upsilon_{d-1}|$, then 
$$
\frac{|\upsilon_{d-1}|}{|\upsilon_{i_0}|}\ge \frac{|\upsilon_{d-1}|}{|\upsilon_{2}|}=\mu
$$
and we are in the case (a) of the step 6. 

If $|\upsilon_2|=|\upsilon_{d-1}|$, then $i_0=1$, we have 
$$
\frac{|\upsilon_{d-1}|}{|\upsilon_1|}\ge \mu
$$
and again we are in the case (a) of the step 6. 

\indent
(b) Assume $|\upsilon_{i_0}|=|\upsilon_d|$. Using Lemma $\ref{Lemma:deuxconjuguesreels}$, we deduce 
$$
|\upsilon_d|>|\upsilon_2|.
$$

If $|\upsilon_2|<|\upsilon_{d-1}|$, then 
$$
\frac{|\upsilon_{i_0}|}{|\upsilon_1|}\ge \frac{|\upsilon_{d-1}|}{|\upsilon_{2}|}=\mu
$$
and we are in the case (b) of the step 6. 

If $|\upsilon_2|=|\upsilon_{d-1}|$, then $i_0=d$, we have 
$$
\frac{|\upsilon_d|}{|\upsilon_2|}\ge \mu
$$
and again we are in the case (b) of the step 6. 

\indent
(c) Assume finally $|\upsilon_1|<|\upsilon_{i_0}|<|\upsilon_d|$. In particular we have $2\le i_0\le d-1$. 
Assume that we are neither in the case (a) nor in the case (b) of the step 6. From 
$$
\frac{|\upsilon_{d-1}|}{|\upsilon_{i_0}|}< \sqrt{\mu}
\quad\hbox{and}\quad
\frac{|\upsilon_{i_0}|}{|\upsilon_2|}< \sqrt{\mu}
$$
we deduce
$$
\frac{|\upsilon_{d-1}|}{|\upsilon_2|}< \mu.
$$
Given the definition of $\mu$, it follows that we have $|\upsilon_2|=|\upsilon_{d-1}|$. Since $\upsilon_{i_0}$ is real, Lemma $\ref{Lemma:deuxconjuguesreels}$ implies $d=3$ and therefore $i_0=2$, $|\upsilon_1|<|\upsilon_2|<|\upsilon_3|$ and
$$
\mu=\min\left\{ 
\frac{|\upsilon_3|}{|\upsilon_2|}, \frac{|\upsilon_2|}{|\upsilon_1|}
\right\}.
$$ 
From
$$
|\gamma_1|\le \lambda_0|\upsilon_1|^a\le \lambda_0\lambda^{-a/2}<1,
\qquad
|\beta_2|=|\beta_{i_0}|<1
$$
and
$$
|\gamma_2 \beta_3|=|\gamma_2 \gamma_3\Psi_3|y\ge \frac{y|\Psi_3|}{|\gamma_1|}\ge \lambda_0^{-1}\lambda^{a/2}|\Psi_3|>1,
$$
we deduce $|\gamma_1 \beta_2|<1<|\gamma_2 \beta_3|$, hence 
$$
\gamma_1\beta_2+ \gamma_2 \beta_3\not=0.
$$
There is an element in the Galois closure of the cubic field  $\Q(\upsilon)$ which maps $\upsilon_1$ to $\upsilon_2$, $\upsilon_2$ to $\upsilon_3$, $\upsilon_3$ to $\upsilon_1$. Therefore, 
$$
\gamma_2\beta_3+ \gamma_3 \beta_1 \not=0.
$$ 
From the part (c) of the step 6 we deduce
$$
0<
\left|
\frac{
\gamma_2\beta_3}{\gamma_3 \beta_1}+1
\right| 
\le 4m \lambda_0^4 \mu^{- a/2}.
$$

\indent
{\tt Step 8}. 
Combining the steps 6 and 7 with the step 4 where we choose
$$
\begin{cases}
i=\ell=d-1, \; j=k=d,\; c=1 &\hbox{in the case (a)},
\\
i=k=i_0, \; j=1,\; \ell=d,\; c=1 &\hbox{in the case (b)},
\\
i=i_0, \; j=k=d,\; \ell=1, \; c=-1 &\hbox{in the case (c)}.
\end{cases}
$$
we deduce
\begin{align}\notag
a\log \mu\le
\Newcst{kappa:step8} R
&
(R+\log |m|+\log \lambda_0+\log \lambda) (\log \lambda) 
\\
\notag
&
\times
\log
\left(
\frac{Ra\log\lambda}{R+\log |m|+\log\lambda_0+\log\lambda} 
\right)\cdotp
\end{align}
For 
$$
 U=\frac{Ra\log\lambda}{R+\log |m|+\log\lambda_0+\log\lambda}
 \quad
 \hbox{and}\quad
 V= \Newcst{kappa:step8b}\frac{R^2(\log\lambda)^2}{\log \mu},
 $$
 we have
 $U\le V\log^\star U$. Therefore we use Lemma $\ref{Lemma:elementary}$ to obtain the conclusion of Theorem \ref{Theorem:majorationdea}.

\section{Proof of Theorem $\ref{Theorem:mainGeneralWeak}$}\label{S:ProofMainTheorem}

Since 
$d\ge 3$, under the assumptions of Lemma $\ref{Lemma:minorationalpha2suralphad}$ we have
$$
\log\frac{|\upsilon''|}{ |\upsilon'|}\ge \frac{ \cst{kappa:lemme1}(\log\lambda)^2
}{ \lambda^{d^2(d+2)/2}}\cdotp
$$
From Lemma $\ref{Lemma:deuxconjuguesreels}$, we deduce that under the assumptions of Theorem $\ref{Theorem:mainGeneralWeak}$ and with the notations of Theorem $\ref{Theorem:majorationdea}$, we have
$$
\log\mu\ge  \frac{  \Newcst{kappa:finale}(\log\lambda)^2
}{ \lambda^{d^2(d+2)/2}}\cdotp  
$$
Hence Theorem $\ref{Theorem:majorationdea}$ implies Theorem $\ref{Theorem:mainGeneralWeak}$.

\bigskip\noindent
{\bf ACKNOWLEDGEMENTS}: This research project was initiated in June 2016 at
the research institute AIMS M'BOUR (Africa) and the authors are grateful to the director,
professor Aissa Wade, for her invitation to give lectures at the center.
The first author benefitted from a grant from NSERC.

\vfill
 \vfill
 
\hbox{
\small
\vbox{
\hbox{Claude \sc Levesque}
	\hbox{D\'{e}partement de math\'{e}matiques et de statistique
	}
	\hbox{Universit\'{e} Laval
	}
	\hbox{Qu\'{e}bec (Qu\'{e}bec)
	}
	\hbox{\sc Canada G1V 0A6
	}
	\hbox{Claude.Levesque@mat.ulaval.ca
	}
}	
\hfill
\vbox{	\hbox{Michel \sc Waldschmidt
	}
	\hbox{Sorbonne Universit\'es
	}
	\hbox{UPMC Univ Paris 06
	}
	\hbox{UMR 7586 IMJ-PRG
	}
	\hbox{F -- 75005 Paris, \sc France 
	}
	\hbox{michel.waldschmidt@imj-prg.fr}
}	
}


\begin{thebibliography}{99}

\providecommand{\bysame}{\leavevmode ---\ }
\providecommand{\og}{``}
\providecommand{\fg}{''}
\providecommand{\smfandname}{\&}
\providecommand{\smfedsname}{\'eds.}
\providecommand{\smfedname}{\'ed.}

\bibitem{BG2debut}
{\sc Y. Bugeaud {\normalfont \smfandname} K. Gy\H{o}ry}, {\em 
Bounds for the solutions of Thue-Mahler equations and norm form equations}, 
\lq\lq Acta Arith.\rq\rq  \,  {\bf 74} (1996), 273--292.
\\
\hbox{ 
\small \tt\url{http://www-irma.u-strasbg.fr/~bugeaud/travaux/BG2debut.ps}}
 
\bibitem{EG}
 {\scshape J-H.~Evertse {\normalfont \smfandname} K.~Gy\H{o}ry} -- {\em
 Unit equation in Diophantine number theory},
 \lq\lq Cambridge studies in advanced mathematics\rq\rq\ \textbf{146}, Cambridge Univ. Press (2015).

\bibitem{GS}
 {\scshape
 X.~Gourdon {\normalfont \smfandname} B.~Salvy} -- {\em Effective asymptotics of linear recurrences with rational coefficients},
 \lq\lq Discrete Mathematics\rq\rq\ \textbf{153} (1996), 145--163.


 \bibitem{LW1}
{\sc C.~Levesque and M.~Waldschmidt}, {\em Familles d'\'equations
 de {T}hue--Mahler n'ayant que des solutions triviales}, \lq\lq Acta Arith.\rq\rq\ {\bf 155} (2012), 117--138.
 \\
\hbox{ 
\small \tt \href{http://arxiv.org/abs/1312.7202}{\tt arXiv:1312.7202 [math.NT]}} 


 \bibitem{LWvdP} 
{\bysame}, 
{\em Families of cubic Thue equations with effective bounds for the solutions},
\lq\lq 
Springer Proceedings in Mathematics \&\ Statistics\rq\rq\ {\bf 43} (2013), 229--243.
 \\
\hbox{ 
\small \tt \href{http://arxiv.org/abs/1312.7204}{\tt arXiv:1312.7204 [math.NT]}}


\bibitem{LW2} 
{\bysame}, 
{\em Solving effectively some families of Thue Diophantine equations}, 
\lq\lq 
Moscow J{.} of Combinatorics and Number Theory\rq\rq\ {\bf 3}, 3--4 (2013), 118--144.
\\
\hbox{ 
\small \tt \href{http://arxiv.org/abs/1312.7205}{\tt arXiv:1312.7205 [math.NT]}}

\bibitem{LW-Ram}
{\bysame}, {\em Familles d'\'equations de Thue associ\'ees \`a 
un sous-groupe de rang $1$ d'unit\'es totalement r\'eelles d'un corps de nombres},
in
\href{http://dx.doi.org/10.1090/conm/655}{SCHOLAR} -- a Scientific Celebration Highlighting Open Lines of Arithmetic Research
(volume dedicated to Ram Murty),
CRM collection
\lq\lq Contemporary Mathematics\rq\rq, 
\href{http://www.ams.org/books/conm/655/}{AMS},
{\bf 655} (2015), 117--134.
\\
\hbox{ 
\small \tt \url{http://www.ams.org/books/conm/655/}}
 \\
\hbox{ 
\small \tt \href{http://arxiv.org/abs/1505.06656}{\tt arXiv: 1505.06656 [math.NT]}}




 \bibitem{LW-Thue}
{\bysame}, {\em
A family of Thue equations involving powers of units of the simplest cubic fields}, 
\lq\lq J. Th\'eor. Nombres Bordx.\rq\rq\ {\bf 27}, No. 2 (2015), 537--563.
 \\
\hbox{ 
\small \tt \href{http://arxiv.org/abs/1505.06708}{\tt arXiv:1505.06708 [math.NT]}}


 

 \bibitem{LW-Balu} 
{\bysame}, 
{\em
Solving simultaneously Thue Diophantine equations: almost totally imaginary case},
\lq\lq Ramanujan Mathematical Society, Lecture Notes Series\rq\rq\
{\bf 23} 
(2016), 137--156.
\\
\hbox{ 
\small \tt \href{http://arxiv.org/abs/1505.06653}{\tt arXiv: 1505.06653 [math.NT]}}




 \bibitem{ST}
{\sc T.N.~Shorey and R.~Tijdeman}, {\em Exponential {D}iophantine equations},
 vol.~87 of Cambridge Tracts in Mathematics, Cambridge University Press,
 Cambridge (1986).

 \bibitem{MR0249395}
{\sc C.L.~Siegel}, {\em
Absch\"atzung von Einheiten,} \lq\lq Nachr. Akad. Wiss. G\"ottingen Math.-Phys.\rq\rq\ Kl. II (1969), 71--86. 
Gesammelte Abhandlungen IV, Springer Collected Works in Mathematics, \S88 (1979), 66--81.

 \bibitem{GL326}
{\sc M.~Waldschmidt}, {\em Diophantine approximation on linear algebraic groups},
 \lq\lq Grundlehren der Mathematischen Wissenschaften\rq\rq\ {\bf 326},
 Springer-Verlag (2000).

\end{thebibliography}
\end{document}